\documentclass[10pt,a4paper]{article}
\usepackage[latin1]{inputenc}
\usepackage{amsmath}
\usepackage{amsfonts}
\usepackage{amssymb}
\usepackage{euler}
 \usepackage[UKenglish]{babel}
\author{Olli Hella}
\title{Algebraic structure of metric value sets}
\date{27 June 2017}
\usepackage{amsthm}
\theoremstyle{plain}
\newtheorem{thm}{Theorem}[section] 
\theoremstyle{definition}
\newtheorem{defn}[thm]{Definition} 
\newtheorem{defns}[thm]{Definitions} 
\newtheorem{exmp}[thm]{Example}
\newtheorem{exmps}[thm]{Examples}
\newtheorem{rmrk}[thm]{Remark}
\newtheorem{lmm}[thm]{Lemma}
\newtheorem{coroll}[thm]{Corollary}
\newtheorem{empt}[thm]{}
\begin{document}
\maketitle
\noindent
ABSTRACT. We introduced the concept of a metric value set (MVS) in an earlier paper \cite{GM}. In this paper we study the algebraic structure of MVSs. For an MVS $M$ we define the concept of $M$-metrizability of a topological space and prove some metrizability results related to those algebraic properties.

\section*{Introduction}
 
First we recall definitions from \cite{GM}. 
A \textit{metric value set} (MVS) is a set $M$ with at least two elements and a binary operation $+$ satisfying the following conditions:
\\
(M1) The operation $+$ is associative.
\\
(M2) The operation $+$ has a (then unique) neutral element $e$. Let $M^{*}=M\setminus\{e\}$.
\\
(M3) If $m_{1}+m_{2}=e$, then $m_{1}=m_{2}=e$.
\\
(M4) For every $m_{1},m_{2}\in M^{*}$ there are $m_{3}\in M^{*}$ and $m_{4},m_{5}\in M$ such that $m_{1}=m_{3}+m_{4}$ and $m_{2}=m_{3}+m_{5}$.

If the operation $+$ is also commutative, we say that $M$ is a \textit{commutative} MVS.

The operation $+$ induces the relations $\unlhd$ and $\lhd$ in $M$ defined by setting $m_{1} \unlhd m_{2} $ if there is $m_{3}\in M$ such that $m_{1}+m_{3}=m_{2}$ and $m_{1}\lhd m_{2}$ if there is $m_{3}\in M^{*}$ such that $m_{1}+m_{3}=m_{2}$.

 A \textit{quasimetric function} $f$ is a map $f\colon X\times X\rightarrow M$ from the square of a set $X$ to an MVS $M$ that satisfies the following conditions:
 \\
(f1) $f(x,z)\unlhd f(x,y)+f(y,z)$ for every $x,y,z\in X$ (triangle inequality).
\\
(f2) $f(x,x)=e$ for every $x\in X$.

A quasimetric function $f$ is a \textit{metric function} if it also satisfies the condition
\\
(f3) $f(x,y)=f(y,x)$ for every $x,y\in X$ (symmetry). 
 
 Let $f\colon X\times X\rightarrow M$ be a quasimetric function. Then the \textit{topology} $\mathcal{T}_{f}$ in $X$ \textit{induced by} $f$ is defined by the family
  \[
B_{f}(x,m)=\{y\in X: f(x,y)\lhd m\},\qquad m\in M^{*},
\]
of \textit{open balls} as the open neighbourhood base at $x$ for every $x\in X$.

In this paper we continue to study the theory of metric value sets. The motivation is to reduce the number of interesting MVSs, i.e., those MVSs that are able to induce different kinds of topologies. The ultimate goal is to be able to characterize whether a given topology can be induced by a quasimetric function. We approach this question by defining morphisms in the category of metric value sets. We define the concept of \textit{$M$-metrizability} to mean that the topology can be induced by a quasimetric function into $M$.

 In this paper we are able to show that every topological space that is $M$-metrizable by a commutative MVS $M$ is $N$-metrizable by a partially ordered MVS $N$. Furthermore we show how to represent MVSs in terms of words and relations. We use the results of \cite{GM} without notification. The study of $M$-metrizability  is continued in the paper \cite{CH}.
 
\section{Metric value sets and partial order}
 
\begin{defn}
Let $(X,\mathcal{T})$ be a topological space. If there exists a quasimetric function $f\colon X\times X \rightarrow M$ 
such that the topology $\mathcal{T}_{f}$ induced by $f$ equals $\mathcal{T}$, then we say that $(X,\mathcal{T})$ or 
shortly $X$ is an $M$-\textit{metrizable} space. We say that $X$ is \textit{MVS-metrizable} if there is an MVS $M$ 
such that $X$ is $M$-metrizable.  
\end{defn} 

\begin{defn}
Let $(M,+)$ be a metric value set and let $\unlhd$ be the induced relation in $M$.
We say that $(M,\unlhd)$ is \textit{antisymmetric} if $\unlhd$ is antisymmetric. Since $\trianglelefteq$ is always reflexive and transitive, we say that $M$ is a \textit{partially ordered set} if it is antisymmetric. If $(M,\unlhd)$ is antisymmetric and for every $m_{1},m_{2}\in M$ we have  $m_{1}\unlhd m_{2}$ or $m_{2}\unlhd m_{1}$, then we say that $M$ is an \textit{ordered set} and $\unlhd$ the \textit{order} of $M$.
\end{defn}

\begin{rmrk}
Even if $M$ is a partially ordered MVS there may still exist elements $m\in M$ for which $m\lhd m$. For example $([ 0,\infty [,\max)$ is an ordered MVS such that $m\lhd m$ for every $m>0$.
\end{rmrk}

\begin{defn} Let $M$ be a set with a binary operation $+$ and $R$ an equivalence relation in $M$. We say that $R$ is a \textit{congruence} in $M$ if from $m R m'$ and $n R n'$ it follows that $(m+n)R(m'+n')$. Let $R$ be a congruence in $M$. We can define a new set $M/R$ whose elements are the  equivalence classes $[m]=\{m'\in M: m'Rm\}$ of the  relation $R$. We say that $M/R$ is the \textit{quotient set} of $M$ by $R$. We define an operation $+_{R}$ in $M/R$ by setting $[m]+_{R}[n]=[m+n]$. This is a well-defined operation because $R$ is a congruence.
\end{defn}

\begin{thm}\label{quotient theorem} \textit{Let $(M,+)$ be an MVS with neutral element $e$ and $R$ a congruence in $M$. Then the following results hold for $(M/R,+_{R})$.} \end{thm}
\noindent
i) \textit{The operation $+_{R}$ is associative. Thus $M/R$ satisfies} (M1).
\\ii) \textit{If $m\unlhd n$, then $[m] \unlhd [n]$.}
\\iii) \textit{The element $[e]$ is a neutral element for $M/R$. Thus $M/R$ satisfies} (M2).
\\iv) \textit{If $[e]=\{e\}$, then $\operatorname{card} M/R\geq 2$ and the pair $(M/R,+_{R})$ satisfies the conditions} (M3) \textit{and} (M4) \textit{and is thus an MVS.}
\\v) \textit{If $(M,+)$ is commutative, then $(M/R,+_{R})$ is commutative.}
\\
\\
\textit{Proof.} i) $([m_{1}]+[m_{2}])+[m_{3}]=[m_{1}+m_{2}]+[m_{3}]=[(m_{1}+m_{2})+m_{3}]=$
\\
$=[m_{1}+(m_{2}+m_{3})]=[m_{1}]+[m_{2}+m_{3}]=[m_{1}]+([m_{2}]+[m_{3}]).$

ii) Let $m+m'=n$. Then $[m]+[m']=[m+m']=[n]$. Thus $[m]\unlhd [n]$.

iii) $[e]+[m]=[e+m]=[m]=[m+e]=[m]+[e]$. Thus $[e]$ is a neutral element for $M/R$.

iv) Let $[e]=\{e\}$. First, there is $m\in M^{*}$, and then $[m]\neq [e]$.

 Now $[m]+[n]=[m+n]=[e]$ if and only if $m+n=e$ which is equivalent to the condition $m=e=n$ and, hence, to the condition $[m]=[e]=[n]$. Thus $M/R$ satisfies the condition (M3).

Let $[m_{1}],[m_{2}]\in (M/R)^{*}$. Then $m_{1},m_{2}\in M^{*}$ and so we find $m_{3}\in M^{*}$ and $m_{4},m_{5}\in M$ for which $m_{3}+m_{4}=m_{1}$ and $m_{3}+m_{5}=m_{2}$. Hence $[m_{3}]\neq [e]$ and $[m_{3}]+[m_{4}]=[m_{1}]$ and $[m_{3}]+[m_{5}]=[m_{2}]$. Therefore (M4) holds for the pair $(M/R,+_{R})$, which is thus an MVS.

v) Trivial to prove. $ \square $

If $M/R$ is an MVS, we say that it is a \textit{quotient metric value set} (QMVS).

\begin{lmm} \textit{If M is a commutative MVS and $m_{1}\unlhd n_{1}, m_{2}\unlhd n_{2},..., m_{k}\unlhd n_{k}$, then $m_{1}+m_{2}+...+m_{k}\unlhd n_{1}+n_{2}+...+n_{k}$. If, in addition, $m_{j}\lhd n_{j}$ for some $j\in \{1,2,...,k\}$, then $ m_{1}+m_{2}+...+m_{k}\lhd n_{1}+n_{2}+...+n_{k}.$} $ \square $
\end{lmm}
\noindent
It is easy to use every MVS $M$ to define the trivial indiscrete topology in any set $X$ by defining a metric function $f\colon X\times X\rightarrow M$ setting $f(x,y)=e$ for every $x,y\in X$. A more interesting question is whether we can use an MVS $M$ to induce a non-trivial topology on some set $X$. The next example shows one way to do it for commutative metric value sets in such way that the image of a quasimetric function $f\colon X\times X \rightarrow M$ is $M$.

\begin{exmp} 
Let $M$ be a commutative MVS. Choose $X=M$. Define $f\colon M\times M \rightarrow M$ by $f(m,n)=e$ when $n\unlhd m$ and $f(m,n)=n$ when $n\ntrianglelefteq m$. Thus $f(n,n)=e$ for every $n\in M$ and so $f$ satisfies the condition (f2).

Let $m_{1},m_{2},m_{3}\in M$. The triangle inequality
\[
f(m_{1},m_{3})\unlhd f(m_{1},m_{2})+f(m_{2},m_{3})
\]
holds trivially if $f(m_{1},m_{3})=e$. Assume $f(m_{1},m_{3})\neq e$. Hence $m_{3}\ntrianglelefteq m_{1}$. Thus  $m_{3}\ntrianglelefteq m_{2}$ or $m_{2}\ntrianglelefteq m_{1}$. If 
$m_{3}\ntrianglelefteq m_{2}$, then $f(m_{2},m_{3})=m_{3}$ and because of the commutativity 
\[
f(m_{1},m_{3})=m_{3}\unlhd f(m_{2},m_{3})+f(m_{1},m_{2})=f(m_{1},m_{2})+f(m_{2},m_{3}).
\]
Assume that $m_{3}\unlhd m_{2}$; then $m_{2}\ntrianglelefteq m_{1}$, whence $f(m_{1},m_{2})+f(m_{2},m_{3})=m_{2}+e=m_{2}\trianglerighteq m_{3}=f(m_{1},m_{3})$. Thus the triangle inequality always holds and $f$ is hence a quasimetric function. Note also that $f(e,m)=m$ for all $m\in M$ and thus $\operatorname{Im}f=M$. 

For a point $m \in M$ the open ball neighbourhoods are of the form
\[
B_{f}(m,n)=\{m'\in M:f(m,m')\lhd n \}\qquad \text{with }n\in M^{*}.
\]
Thus $m'\in B_{f}(m,n)$ if and only if $m'\unlhd m$ or $m'\lhd n$. We see that if $m'\unlhd m$, then $m'$ belongs to every neighbourhood of $m$.

For example if $(M,+)=([0,\infty[,+)$, then $f$ induces the topology
\[
\mathcal{T}_{f}=\{A\subset [0,\infty[ : A \text{ is an interval and } 0\in A \neq \{0\}\}.
\]
\end{exmp}

\begin{empt} \textbf{Constructing a partially ordered MVS from a commutative MVS.}\label{constructing partially}
Let $M$ be a commutative MVS. Define an equivalence relation $R$ in $M$ by setting $mRn$ if and only if $m\unlhd n$ and $n\unlhd m$. The relation $R$ is clearly symmetric, and its reflexivity and  transitivity follow from the corresponding attributes of $\unlhd$.

The relation $R$ is compatible with $+$: If $m,n,m',n'\in M$, $mRm'$ and $nRn'$, then
\[
(m+m'')+(n+n'')=m'+n'
\]
with some $m'',n''\in M$ and from the commutativity of $M$ we see that $m+n\unlhd m'+n'$. Symmetrically we see that $m'+n'\unlhd m+n$ so $(m+n)R(m'+n')$. Thus we get a well-defined operation $+_{R}$ to the quotient set $M/R$ by setting $[m]+_{R}[n]=[m+n]$. Since $[e]=\{e\}$, the pair $(M/R,+_{R})$ is an MVS by \ref{quotient theorem}iv). Denote the induced relations in $M/R$ by $\lhd_{R}$ and $\unlhd_{R}$.
\\
\\
\underline{Claim}. \textit{If $[m]\unlhd_{R} [n]$, then $m\unlhd n$. If $[m]\lhd_{R} [n]$, then $m\lhd n$.} 
\\
\\
\textit{Proof.} Let $[m]\unlhd_{R} [n]$. Then there exists $m'\in M$ such that $[m]+_{R}[m']=[n]$. Therefore $[m+m']=[n]$. Thus $(m+m')Rn$ and so $m+m'\unlhd n$. Thus $m\unlhd n.$ If $[m]\lhd_{R}[n]$, then we may assume that $m'\in M^{*}$, and hence $m\lhd n$ $\square$
\\
\\
Let $[m],[n]\in M/R$, $[m]\unlhd_{R}[n]$ and $[n]\unlhd_{R}[m]$. Now $mRn$ and therefore $[m]=[n]$. Thus $\unlhd_{R}$ is antisymmetric and therefore it is a partial order. Hence  $(M/R,+_{R})$ is a partially ordered MVS constructed from the commutative MVS $M$. 

Later in \ref{partial theorem} we show that when $M$ is a commutative MVS, then every $M$-metrizable space is $M/R$-metrizable. Thus in the category of commutative MVSs we can focus our study to partially ordered MVSs.
\end{empt}

\begin{defn}
Let $(M,+)$ be an MVS and $N$ a subset of $M$ that contains $e$ and is closed with respect to the operation $+$. Let $\operatorname{card}N\geq 2$ and let the restriction $+|_{N\times N}\colon N\times N \rightarrow N$ of $+$ to the set $N\times N$ satisfy the condition (M4). Since $(N, +|_{N\times N})$ satisfies the conditions (M1)--(M3) trivially, it is an MVS. We say that $(N, +|_{N\times N})$ or shorter $N$ is a \textit{sub-metric value set} or a \textit{sub-MVS} of $M$. From the definition we easily see that if $N$ is a sub-MVS of $M$ and $O$ is a sub-MVS of $N$, then $O$ is a sub-MVS of $M$.
\end{defn}

\section{Morphisms in metric value sets }
In this section we study maps between metrizable spaces and maps between metric value sets.

\begin{defn}[\textbf{Equivalence of quasimetric functions}] Let $X$ be a set, $M_{1},M_{2}$ metric 
value sets and $f_{i}\colon X\times X \rightarrow M_{i}$ a quasimetric function for $i=1,2$. 
We say that the quasimetric function $f_{2}$ is \textit{finer} than $f_{1}$ if $\mathcal{T}_{f_{1}}\subset \mathcal{T}_{f_{2}}$. This is true if and only if for every $x\in X$ it holds that for every $m_{1}\in M_{1}^{*}$ there exists $m_{2}\in M_{2}^{*}$ such that $B_{f_{2}}(x,m_{2})
\subset B_{f_{1}}(x,m_{1})$, i.e., such that for every $y\in X$ we have
\[
f_{2}(x,y)\lhd m_{2}\Longrightarrow f_{1}(x,y)\lhd m_{1}.
\] 
We say that the quasimetric functions $f_{1}$ and $f_{2}$ are \textit{equivalent} if $\mathcal{T}_{f_{1}}=\mathcal{T}_{f_{2}}$. Thus $f_{1}$ and $f_{2}$ are equivalent if and only if $f_{1}$ is finer than $f_{2}$ and $f_{2}$ is finer than $f_{1}$.

\end{defn}
\begin{defn}[\textbf{MVS homomorphisms}] Let $(M_{1},+_{1})$ and $(M_{2},+_{2})$ be two metric value sets and $h\colon M_{1}\rightarrow M_{2}$ a map. We say that $h$ is a \textit{metric value set homomorphism} if it satisfies the following two conditions:
\\
(H1) For all $m_{1}\in M$ it holds that $h(m_{1})=e_{M_{2}} \Longleftrightarrow m_{1}=e_{M_{1}}$.
\\
(H2) $h(m+_{1}n)=h(m)+_{2}h(n)$ for every $m,n\in M_{1}$.
\\
Since there is no danger of misunderstanding we may shortly say that $h$ is a homomorphism.
\end{defn}

\begin{exmps}\label{homomorphism examples} 
1. The map $h\colon (\mathbb{N},+)\rightarrow (\{0,1\},\max)$, $0\mapsto 0$, $n\mapsto 1$ for $n\neq 0$, is a homomorphism

2. The map
 $\operatorname{Id}_{M}\colon M\rightarrow M$, $m\mapsto m$, is a homomorphism for every MVS $M$.
 
 3. Let $N$ be a sub-MVS of an MVS $M$. Then the inclusion $i\colon N \hookrightarrow M$, $n\mapsto n$, is a homomorphism.
 
 4. Let $M$ be an MVS and $R$ a congruence in $M$. Then $M/R$ is a QMVS and $h\colon M\rightarrow M/R$, $m\mapsto [m]_{R}$, a homomorphism if and only if $[e_{M}]=\{e_{M}\}$.
\end{exmps}

\begin{lmm}\label{homomorphism lemma} \textit{Let $h\colon M_{1}\rightarrow M_{2}$ be a homomorphism and $m,n\in M_{1}$. If $m\lhd n$, then $h(m)\lhd h(n)$, and if $m\unlhd n$, then $h(m)\unlhd h(n)$.} $ \square $
\end{lmm}
\begin{lmm} \textit{Let $h\colon M_{1}\rightarrow M_{2}$ and $g\colon M_{2}\rightarrow M_{3}$ be homomorphisms. Then $g\circ h\colon M_{1}\rightarrow M_{3}$ is a homomorphism.}
\end{lmm}
\noindent
\textit{Proof.} First we notice that
\[
(g\circ h)(e_{M_{1}})=g(h(e_{M_{1}}))= g(e_{M_{2}})=e_{M_{3}}
\]
and when $m_{1}\neq e_{M_{1}}$, then $h(m_{1})\neq e_{M_{2}}$ and thus
\[
(g\circ h)(m_{1})=g(h(m_{1}))\neq e_{M_{3}}.
\]
Furthermore for all $m,n\in M_{1}$ we have
\begin{align*}
(g\circ h)(m+n)=& g(h(m+n))=g(h(m)+h(n))
\\
=& g(h(m))+g(h(n))=(g\circ h)(m)+ (g\circ h)(n). \quad \square
\end{align*}

\begin{thm} Let $h\colon M\rightarrow N$ be a homomorphism. Then the image $h(M)$ is a sub-MVS of $N$.
\end{thm}
\noindent
\textit{Proof.} Let $m_{i}\in M$ and $n_{i}=h(m_{i})$ for $i=1,2$. Then $n_{1}+n_{2}=h(m_{1})+h(m_{2})=h(m_{1}+m_{2})$. Therefore $h(M)+h(M)\subset h(M)$. 

We use the symbol $+$ also to mean the restriction of $+$ to $h(M)$. It is associative as a restriction of an associative operation. The image of $e_{M}$ under $h$ is $e_{N}$, the neutral element of $h(M)$. The pair $(h(M),+)$ also satisfies the condition (M3) and the condition $\operatorname{card}h(M)\geq 2$, because $h^{-1}\{e_{N}\}=\{e_{M}\}$. The condition (M4) holds because for every $h(m_{1}),h(m_{2})\in h(M)^{*}$ we can choose an appropriate $h(m_{3})\in h(M)^{*}$ by requiring that $e_{M}\neq m_{3}\unlhd m_{1},m_{2}$. Hence $h(M)$ is a sub-MVS of $N$. $ \square $ 

\begin{defn} Let $h\colon M\rightarrow N$ be a homomorphism. Define a relation $\ker(h)$ in $M$ by setting
\[
(m_{1},m_{2})\in \ker(h) \Longleftrightarrow h(m_{1})=h(m_{2}).
\]
It is easy to check that $\ker(h)$ is an equivalence relation in $M$. The relation $\ker(h)$ is also compatible with $+$: Let $(m_{1},m_{1}'),(m_{2},m_{2}')\in \ker(h)$. Then
\[
h(m_{1}+m_{2})=h(m_{1})+h(m_{2})= h(m_{1}')+h(m_{2}')=h(m_{1}'+m_{2}').
\]
Thus $(m_{1}+m_{2},m_{1}'+m_{2}')\in \ker(h)$; so $\ker(h)$ is a congruence. We say that $\ker(h)$ is the \textit{kernel} of the homomorphism $h$. The set $M/ker(h)$ is a QMVS by \ref{quotient theorem}iv) because from $m \ker(h) e_{M}$ it follows that $h(m)=h(e_{M})=e_{N}$ and thus $m=e_{M}$. It also follows that the map $q:M\rightarrow M/ker(h)$, $m\mapsto [m]$, is a homomorphism. 
\end{defn}

\begin{thm}\label{metric transform} Let $M_{1}$ and $M_{2}$ be two metric value sets and $f\colon X\times X \rightarrow M_{1}$ a quasimetric function. Let $h\colon M_{1}\rightarrow M_{2}$ be a homomorphism. Then $h\circ f\colon X\times X\rightarrow M_{2}$ is a quasimetric function. If $f$ is a metric function, then $h\circ f$ is a metric function
\end{thm}
\noindent
\textit{Proof.}
The condition (f2) holds because $(h\circ f)(x,x)=h(f(x,x))=h(e_{M_{1}})=e_{M_{2}}$. Let $x,y,z\in X$. Then $f(x,z)\unlhd f(x,y)+f(y,z)$ and hence
\begin{align*}
(h\circ f)(x,z)=& h(f(x,z)) \unlhd h(f(x,y)+f(y,z)) 
\\
=& h(f(x,y))+h(f(y,z))=(h\circ f)(x,y)+(h\circ f)(y,z).
\end{align*}
Thus (f1) holds. 

Assume that $f$ is symmetric. Then $h\circ f$ is symmetric. Thus if $f$ is a metric function, then also $h\circ f$ is such. $\square$
\\
\\
It follows that also $h\circ f$ induces a topology $\mathcal{T}_{h\circ f}$ in $X$. This topology is not necessarily the same as $\mathcal{T}_{f}$. However $\mathcal{T}_{h\circ f}$ is always the same as $\mathcal{T}_{f}$ if $h$ satisfies certain conditions not depending on $f$. If there is $h\colon M_{1}\rightarrow M_{2}$ satisfying those conditions, then every $M_{1}$-metrizable space 
is $M_{2}$-metrizable. Next we examine what those conditions are. 

\begin{defn} Let $h\colon M_{1} \rightarrow M_{2}$. We say that $h$ is an \textit{isomorphism} of metric value sets if it is a bijective homomorphism.
\end{defn}

\begin{thm}The inverse of an isomorphism is an isomorphism.
\end{thm}
\noindent
\textit{Proof.} Let $h\colon M_{1}\rightarrow M_{2}$ be an isomorphism. The inverse $h^{-1}\colon M_{2}\rightarrow M_{1}$ is a bijection, so it is enough to show it is a homomorphism. Clearly $h^{-1}(e_{M_{2}})=e_{M_{1}}$. Let $m,n\in M_{2}$. Now
\begin{align*}
h^{-1}(m)+h^{-1}(n)=& h^{-1}(h(h^{-1}(m)+h^{-1}(n)))
\\
=& h^{-1}(h(h^{-1}(m))+h(h^{-1}(n)))=h^{-1}(m+n).
\end{align*}
Therefore $h^{-1}$ is a homomorphism, so it is an isomorphism. $ \square $
\\
\\
If there exists an isomorphism from $M$ onto $N$, then we say that those MVSs are \textit{isomorphic} and we denote $M\cong N$.  

\begin{thm} Let $h\colon M\rightarrow N$ be a homomorphism. Then $M/\ker(h)$ and $h(M)$ are isomorphic.
\end{thm}
\noindent
\textit{Proof.} Define a map $h^{*}\colon M/\ker(h)\rightarrow h(M)$ by setting $h^{*}([m])=h(m)$, where $[m]$ is the equivalence class of an element $m\in M$ with respect to $\ker(h)$. The map $h^{*}$ is well-defined. It is a homomorphism, because $h^{*}([m])=h(m)=e_{N}$ if and only if $m=e_{M}$, i.e., $[m]=[e_{M}]$, and
\[
h^{*}([m])+h^{*}([n])=h(m)+h(n)=h(m+n)=h^{*}([m+n])=h^{*}([m]+[n]).
\]
It is an injection, because the equivalence classes are formed exactly from those elements whose image in $h$ is the same element. It is a surjection because every $m\in M$ belongs to some equivalence class. Thus $h^{*}$ is the isomorphism that we were looking for. $\square$

\begin{rmrk} Homomorphisms, isomorphisms, sub-MVSs and QMVSs are named as in the universal algebra. Still we can not directly use results of universal algebra to MVSs because the conditions (M3) and (M4) can not be stated without using universal and existential quantifiers. Nonetheless many algebraic properties of MVSs can be proven in a similar way as the corresponding results in universal algebra \cite{UA}.
\end{rmrk}

\begin{defn} We say that a homomorphism $h\colon M_{1}\rightarrow M_{2}$ is \textit{fine}, if for every $m_{2}\in M_{2}^{*}$ there exists $m_{1}\in M_{1}^{*}$ such that $h(m_{1})\unlhd m_{2}$.
\end{defn}

\begin{lmm} \textit{An isomorphism is fine.}
\end{lmm}
\noindent 
\textit{Proof.} Let $h\colon M_{1}\rightarrow M_{2}$ be an isomorphism. Let $m_{2}\in M_{2}^{*}$. Since $h$ is a bijection, we see that there exists $h^{-1}(m_{2})\in M_{1}^{*}$, and $h(h^{-1}(m_{2}))=m_{2}\unlhd m_{2}$. Therefore $h$ is fine. $ \square $

For the next theorem, recall the MVS $M_{\infty}$ in [2, 1.3.8].

\begin{thm}\label{infty theorem} Let $M$ be an MVS and $i\colon M\hookrightarrow M_{\infty}$ the inclusion. Then $i$ is a fine homomorphism.
\end{thm}
\noindent
\textit{Proof.} Notice that $M$ is a sub-MVS of $M_{\infty}$. By \ref{homomorphism examples}.3, $i$ is a homomorphism. Let $m\in M^{*}_{\infty}.$ If $m\neq \infty$, then $i(m)=m\unlhd m$. Let $m=\infty$. Choose $n\in M^{*}$; then $i(n)+_{\infty}\infty=\infty.$ Therefore $i(n)\unlhd \infty$. Thus $i$ is fine. $\square$

\begin{thm}\label{inclusion theorem} If $f\colon X\times X\rightarrow M_{1}$ is a quasimetric function and $h\colon M_{1} \rightarrow M_{2}$ a fine homomorphism, then $\mathcal{T}_{h\circ f} \subset \mathcal{T}_{f}$. 
\end{thm}
\noindent
\textit{Proof.} Let $x\in X$ and $m_{2}\in M_{2}^{*}$. Let $m_{1}\in M_{1}^{*}$ be such that $h(m_{1})\unlhd m_{2}$. The following deduction shows that $B_{f}(x,m_{1})\subset B_{h\circ f}(x,m_{2})$:
\begin{align*}
y\in B_{f}(x,m_{1}) &\Longrightarrow f(x,y)\lhd m_{1} 
\\
&\Longrightarrow(h\circ f)(x,y)\lhd h(m_{1})\unlhd m_{2}
\Longrightarrow y\in B_{h\circ f}(x,m_{2}).
\end{align*}

From this it follows that $\mathcal{T}_{h\circ f}\subset \mathcal{T}_{f}$. $\square$

\begin{coroll} \textit{Let $f\colon X\times X \rightarrow M_{1}$ be a quasimetric function and $h\colon M_{1} \rightarrow M_{2}$ an isomorphism. Then $\mathcal{T}_{f} = \mathcal{T}_{h\circ f}$}. $\square$
\end{coroll}

\begin{thm} Let $h\colon M_{1} \rightarrow M_{2}$ be an isomorphism. Then a topological space is $M_{1}$-metrizable if and only if it is $M_{2}$-metrizable.
\end{thm}
\noindent
\textit{Proof.} Let $f\colon X\times X \rightarrow M_{1}$ be a quasimetric function. Then $h\circ f\colon X \times X \rightarrow M_{2}$ is a quasimetric function by \ref{metric transform} and furthermore $\mathcal{T}_{f}=\mathcal{T}_{h\circ f}$. Therefore if $(X,\mathcal{T})$ is $M_{1}$-metrizable, then $(X,\mathcal{T} )$ is $M_{2}$-metrizable. Every $M_{2}$-metrizable space is $M_{1}$-metrizable, because $h^{-1}\colon M_{2}\rightarrow M_{1}$ is also an isomorphism. $\square$

\begin{exmps} 1. The map 
\[
h\colon([0,\infty[,+) \rightarrow ([0,\infty[,+),\qquad m\mapsto 2m ,
\]
is an isomorphism.
\\
2. The map
\[
h\colon([0,\infty[,+) \rightarrow ([1,\infty[,\cdot),\qquad m\mapsto e^{m},
\]
is an isomorphism. Therefore every $([0,\infty[,+)$-metrizable space is also $([1,\infty[,\cdot)$-metrizable, and conversely.
\end{exmps}

\begin{thm} Every $M$-metrizable space, where $M$ is commutative, is also $N$-metrizable with some partially ordered commutative metric value set $N$.\label{partial theorem}
\end{thm}
\noindent
\textit{Proof.} Let $M$ be a commutative MVS and $X$ an $M$-metrizable topological space. Let the topology of $X$ be induced by a quasimetric function $f\colon X\times X\rightarrow M$. Let $M/R$ be as in \ref{constructing partially}. Then the map $h\colon M\rightarrow M/R, m\mapsto [m]$, is a homomorphism. The homomorphism $h$ is fine because for every $[m]\in 
(M/R)^{*}$ it holds that $m\in M^{*}$ and $h(m)=[m]\unlhd [m]$. Therefore $\mathcal{T}_{h\circ f}\subset \mathcal{T}_{f}$

Let $x\in X$ and $m\in M^{*}$. Let $y\in B_{h\circ f}(x,[m])$. Denote $[m'] = (h\circ f)(x,y)\lhd [m]$. Then $f(x,y)Rm'$, so $f(x,y)\unlhd m'$, and on the other hand $m' \lhd m$ by the Claim in \ref{constructing partially}. Thus $f(x,y)\lhd m.$ Therefore $y \in B_{f}(x,m)$ proving that $B_{h\circ f}(x,[m])\subset B_{f}(x,m)$. Thus $\mathcal{T}_{f}\subset \mathcal{T}_{h\circ f}$. Hence $\mathcal{T}_{f}=\mathcal{T}_{h\circ f}$. $\square$

\begin{thm} Let $f\colon X\times X\rightarrow M$ be a quasimetric function and $i\colon M\rightarrow M_{\infty}$ the inclusion. Then $\mathcal{T}_{f}=\mathcal{T}_{i\circ f}$.
\end{thm}
\noindent
\textit{Proof.} By \ref{infty theorem} and \ref{inclusion theorem} we have $\mathcal{T}_{i\circ f}\subset \mathcal{T}_{f}$. Let $x\in X$ and $m\in M^{*}$. Then $i(m)\in M_{\infty}^{*}$ and $B_{i\circ f}(x,i(m))=B_{f}(x,m)$. Therefore $\mathcal{T}_{f}\subset \mathcal{T}_{i\circ f}$. Thus $\mathcal{T}_{f}=\mathcal{T}_{i\circ f}$. $\square$

\begin{defns} To construct a representation theory for MVSs we introduce the following two concepts. Let $V$ be a non-empty set. A \textit{word} is a finite sequence $(v_{1},...,v_{n})$ of elements of $V$. Denote shortly $v_{1}...v_{n}=(v_{1},...,v_{n})$. An empty sequence, with $n=0$, is called an \textit{empty word} and denoted as $0$. We denote the collection of all words of the set $V$ with the symbol $\textbf{W}(V)$. We define an operation $*$ in $\textbf{W}(V)$ called \textit{concatenation} in the following way. For every two words $v=v_{1}...v_{n}\in \textbf{W}(V)$ and $v'=v_{1}'...v_{k}'\in \textbf{W}(V)$ the concatenation $v*v'$ is $v_{1}...v_{n}v_{1}'...v_{k}'\in \textbf{W}(V)$.

We also define the \textit{transitive closure}. Let $X$ be a set and $R$ a relation in $X$. The transitive closure $R_{T}$ is the smallest transitive relation containing $R$. It is formed in the following way. Set $R_{0}=R$ and
\[
R_{i+1}=R_{i} \cup \{(a,c):\exists b\in X \text{ such that } (a,b),(b,c)\in R_{i} \}.
\]
Now $R_{T}$ is $\bigcup_{i\geq 0}R_{i}$. 
\end{defns}

\begin{empt}\textbf{Relations in the set of words.}
Let $V$ be a non-empty set. Concatenation $*$ in $\textbf{W}(V)$ satisfies the conditions (M1), (M2) and (M3), but not the condition (M4), unless $V$ has only one element.

Let $a,b,c\in V$. Let $R_{ab\sim c}$ be the relation in $\textbf{W}(V)$ for which $vR_{ab\sim c}u$ if and only if $u=v$ or the word $u$ can be formed from $v$ by replacing $ab$ in $v$ in one place with $c$ or $c$ in one place with $ab$. More precisely $uR_{ab\sim c}v$ if and only if $v=u$ or $v$ and $u$ can be presented in the form $v=v_{1}...v_{k}abv_{k+1}...v_{n}$ and $u=v_{1}...v_{k}cv_{k+1}...v_{n}$ or in the form $v=v_{1}...v_{k}cv_{k+1}...v_{n}$ and $u=v_{1}...v_{k}abv_{k+1}...v_{n}$.

The relation $R_{ab\sim c}$ is symmetric and reflexive, but not necessarily transitive. Let $\textbf{R}$ be a non-empty union of a family of the relations $R_{ab\sim c}$ and of the relation $\operatorname{id}_{\textbf{W}(V)}=\{(a,a): a\in \textbf{W}(V)\}$. Thus $\textbf{R}$ is symmetric and reflexive. Hence the transitive closure $\textbf{R}_{T}$ is also symmetric and reflexive and therefore an equivalence relation in $\textbf{W}(V)$.
\end{empt}

\begin{exmp} Let $V=\{a,b,c,d\}$ and $\textbf{R}=R_{ab\sim c}\cup R_{ad\sim b} \cup R_{bc\sim a}$. Then $(abcd,c)\in \textbf{R}_{T}$, but $(abcd,d)\notin \textbf{R}_{T}$. The former is true due to the following relation chain:
\[
abcdR_{bc\sim a}aadR_{ad\sim b}abR_{ab\sim c}c.
\]
The latter holds because $d$ is in relation $R$ only with itself for each $R\in\{R_{ab\sim c},
\\
 R_{ad\sim b}, R_{bc\sim a}\}$
\end{exmp}

\begin{thm}The relation $\textbf{R}_{T}$ is a congruence in the set $\textbf{W}(V)$.
\end{thm}
\noindent
\textit{Proof.}
Let $u\textbf{R}_{T}v$ and $u'\textbf{R}_{T}v'$. Then $uu'\textbf{R}_{T}vu'$ and $vu'\textbf{R}_{T}vv'$. Thus $uu'\textbf{R}_{T}vv'$. $ \square $ 
\begin{thm}\label{M4 theorem} The pair $(\textbf{W}(V)/\textbf{R}_{T},*)$ is an MVS exactly when for every $u,v\in V$ there exist $w\in V$ and $\bar{u},\bar{v}\in \textbf{W}(V)$ such that $w\bar{u}\textbf{R}_{T}u$ and $w\bar{v}\textbf{R}_{T}v$.
\end{thm}
\noindent
\textit{Proof.} The pair $(\textbf{W}(V)/\textbf{R}_{T},*)$ satisfies the conditions (M1) and (M2), because $\textbf{W}(V)$ satisfies them and $\textbf{R}_{T}$ is a congruence. The condition (M3) holds 
since $[0]=\{0\}$. For the same reason we see that $\textbf{W}(V)/\textbf{R}_{T}$ contains at least two elements.

Assume the condition. Let $[a],[b]\in \textbf{W}(V)/\textbf{R}_{T}\setminus \{[0]\}$. Choose $u,v\in V$ and $\bar{a},\bar{b}\in \textbf{W}(V)$ with $a=u\bar{a}$ and $b=v\bar{b}$.
The condition yields $w\in V$ and $\bar{u},\bar{v}\in \textbf{W}(V)$ such that $[u]=[w][\bar{u}]$ and $[v]=[w][\bar{v}]$. Therefore $[w]\neq [0], [a]=[w][\bar{u}\bar{a}]$ and $[b]=[w][\bar{v}\bar{b}]$. Thus (M4) holds, and hence $(\textbf{W}(V)/\textbf{R}_{T},*)$ is an MVS.

Assume that $(\textbf{W}(V)/\textbf{R}_{T},*)$ is an MVS. Let $u,v\in V$. By (M4) choose $a,u',v'\in \textbf{W}(V)$ with $[a]\neq [0]$, $[u]=[a][u']$ and $[v]=[a][v']$. Choose $w\in V$ and $\bar{a}\in \textbf{W}(V)$ with $a=w\bar{a}$. Let $\bar{u}=\bar{a}u'$ and $\bar{v}=\bar{a}v'$. Then $[u]=[w][\bar{u}]$ and $[v]=[w][\bar{v}]$. Thus the condition holds. $\square$

\begin{rmrk} On the basis of the previous theorem it is enough to study whether for every $a,b\in V,$ $a\neq b$, there exists  $c\in V$ such that $[c]\unlhd [a],[b]$. If $a=b$, then $[a]\unlhd [a]=[b]$.
\end{rmrk}

\begin{exmps} 1. Let $V=\{a\}$ and $\textbf{R}=\text{id}_{\textbf{W}(V)}$. Now $\textbf{W}(V)/\textbf{R}_{T}$ can be identified with $\textbf{W}(V)$ and the elements of $\textbf{W}(V)$ are the finite sequences $a...a$ and the empty word $0$.

2. Let $V=\{a,b\}$ and $\textbf{R}=\text{id}_{\textbf{W}(V)}$. Now $(\textbf{W}(V)/R_{T},*)$ does not satisfy the condition in the previous theorem, because $[a]\ntrianglelefteq[b]$ and $[b]\ntrianglelefteq[a]$.

3. Let $V=\{a,b\}$ and $\textbf{R}=R_{ab\sim b}$. Now $[a]\unlhd [a]$ and also $[a]\unlhd [ab]=[b]$. Therefore $(\textbf{W}(V)/R_{T},*)$ is an MVS. Its elements are
\[
[0],[a],[b],[aa],[ba],[bb],[aaa],[baa],[bba],[bbb],[aaaa],[baaa],...\,.
\] 

4. Let $X$ be a set and $V=\{1\}\cup X$ with $1\notin X$. Let $\textbf{R} =(\bigcup_{x\in X}R_{1x\sim x})\cup \operatorname{id}_{\textbf{W}(V)}$. Now $[1]\unlhd [x]$ for all $x\in X$. Therefore  $(\textbf{W}(V)/\textbf{R}_{T},*)$ is an MVS.

5. Let $V=\{1,a,b,c,d\}$ and
\[
\textbf{R}=R_{1a\sim a}\cup R_{1b\sim b}\cup R_{1c\sim c}\cup R_{1d\sim d}\cup R_{ab\sim c}\cup R_{ba\sim d}.
\]
The pair $(\textbf{W}(V)/\textbf{R}_{T},*)$ is an MVS, because $[1]\unlhd [x]$ for every $x\in V$. 

Define a map $I\colon\textbf{W}(V)\rightarrow \{-1,0,1\}$ by $I(v)=0$ if none of the letters $a,b,c$ and d occurs in $v$ and otherwise by $I(v)=1$ or $I(v)=-1$, respectively, if the first of them occurring in $v$ is $a$ or $c$ or, respectively, $b$ or $d$.

We see that $u\textbf{R}v\Longrightarrow I(u)=I(v).$
Therefore $[a][b]=[ab]\neq [ba]=[b][a]$, because $I(ab)=1\neq -1=I(ba)$. Thus $\textbf{W}(V)/\textbf{R}_{T}$ is not commutative.
\end{exmps}

\begin{thm}[\textbf{Representation theorem for MVSs}]
\textit{Every MVS is isomorphic to an MVS of the form $(\textbf{W}(V)/\textbf{R}_{T},*)$.}
\end{thm}
\noindent
\textit{Proof.} Let $(M,+)$ be an MVS. Choose $V=M^{*}$ and 
\[
\textbf{R}= \bigcup\{ R_{m_{1}m_{2}\sim m_{3}}:m_{1},m_{2},m_{3}\in M^{*}, m_{1}+m_{2}=m_{3} \}.
\]
We next show that the pair $(\textbf{W}(V)/\textbf{R}_{T},*)$ is an MVS by checking the condition in \ref{M4 theorem}. Thus, let $m_{1},m_{2}\in M^{*}$. By (M4) choose $m_{3}\in M^{*}$ and $m_{4},m_{5}\in M$ with $m_{1}=m_{3}+m_{4}$ and $m_{2}=m_{3}+m_{5}$. Then $[m_{1}]=[m_{3}m_{4}]=[m_{3}][m_{4}]$ and $[m_{2}]=[m_{3}m_{5}]=[m_{3}][m_{5}]$.

Define a map $h\colon M\rightarrow \textbf{W}(V)/\textbf{R}_{T}$ by $m\mapsto [m]$ for every $m\in M^{*}$ and by $e_{M}\mapsto [0]$. It is a surjection, because every word of length $n\geq 2$ is in relation $\textbf{R}$ with some word with length $n-1$. Thus by induction we see that every word is in relation $\textbf{R}_{T}$ with a word of length one or zero.

The map $h$ is an injection: Let $m_{1},m_{2}\in M^{*}$ and $[m_{1}]=[m_{2}]$. Let 
\[
m_{1}Rn_{1,1}...n_{1,k_{1}}Rn_{2,1}...n_{2,k_{2}}R...Rm_{2}
\]
be the corresponding relation chain. Now
\[
m_{1}=n_{1,1}+...+n_{1,k_{1}}=n_{2,1}+...+n_{2,k_{2}}=...=m_{2}.
\]
Therefore $m_{1}=m_{2}$. Furthermore $h(m)=[0]$ if and only if $m=e_{M}$. Thus $h$ is an injection.

The map $h$ is a homomorphism, because in addition $h(m_{1})*h(m_{2})=[m_{1}]*[m_{2}]=[m_{1}m_{2}]=h(m_{1}+m_{2})$ if $m_{1},m_{2}\in M^{*}$, $h(m)*h(e_{M})=[m]*[0]=[m0]=[m]=h(m)=h(m+e_{M})$ and similarly $h(e_{M})*h(m)=h(e_{M}+m)$ if $m\in M^{*}$, and finally $h(e_{M})*h(e_{M})=[0]*[0]=[00]=[0]=h(e_{M})=h(e_{M}+e_{M})$. Therefore $h$ is an isomorphism. $\square$

\begin{empt}\textbf{A question for further study.}
It is an open question whether there exists a topological space $X$ such that $X$ is $M$-metrizable by a non-commutative MVS $M$ but not $N$-metrizable by any commutative MVS $N$.

\end{empt}

\section*{Acknowledgement}
The author wants to thank Jouni Luukkainen for many comments and suggestions that helped to improve the structure,  clarity of proofs and language in this paper.

\end{document}